# A Robust Power Grid Defense Model Considering Load Demand and Wind Generation Uncertainties

Yingmeng Xiang*,§, Xiaohu Zhang*, Di Shi*, Yanming Jin†, Zhiwei Wang* and Lingfeng Wang§

*GEIRI North America, San Jose, CA 95134, USA
§University of Wisconsin, Milwaukee, WI 53211, USA
†State Grid Energy Research Institute, Beijing, China

*Abstract*—It is a major task to develop effective strategies for defending the power system against deliberate attacks. It is critical to comprehensively consider the human-related and environmental risks and uncertainties, which is missing in existing literature. This paper considers the load demand uncertainties and wind generation uncertainties in addition to the interactive attacker/defender behaviors. Specifically, a defender-attacker-nature-operator model is proposed, which incorporates the attack/defense interaction, the corrective re-dispatch of the operator, the coordination between the attack strategy and the stochastic nature of load demands and wind generations. The Column-and-Constraint Generation (C&CG) algorithm is adopted for solving the proposed model by decomposing the proposed model into a master problem and a sub-problem. Simulations are performed using MATLAB and CPLEX on a modified IEEE RTS79 system. The simulation results verify the validity of the proposed model.

I. INTRODUCTION

The cyber-physical security of power systems is a major concern in operation and planning, as the power grid plays a critical role in supporting our modern human society. In recent years, with the adoption of the cutting-edge smart grid technologies, such as renewable generation, microgrid, various kinds of flexible alternating current transmission systems (FACTS) devices, advanced metering systems [1]-[4], etc., the power system is rapidly evolving into a complicated interconnected cyber-physical system involving human in the loop. They can greatly improve the power system economy and flexibility, yet inevitably bring more uncertainties and risks of attacks to the power system, which may weaken the power system's secure and reliable operation.

In the past years, several severe attacks have launched against the power grid, which caused tremendous losses. Two typical examples were the cyber attack on the Ukrainian power grid in 2015 [5] and the shooting at the substation in Silicon Valley, USA in 2013 [6]. These kinds of attacks may be more frequent and disastrous in the future. Thus, it is urgent to defend the power system against man-made attacks. Although the N-1 or even N-2 reliability criterion is enforced in the power grid, it is insufficient to deal with attacks. Such criteria were adopted mainly to deal with the random failures of major power system components, but not capable of handling man-made attacks which may cause the simultaneous outages of multiple elements.

When an attack is imminent, it is usually not affordable to defend all the components in a bulk power grid; therefore, it is important to identify the key elements. Ideally, the power system operator needs to consider the risks/uncertainties arising from the nature (e.g., weather, temperature), the stochastic behaviors of customers, attacks from the cyber domain and physical domain, the uncertain characteristics of the power system devices, etc. In the existing literature, some researchers proposed interesting models and ideas regarding the identification of critical components. The bilevel attacker-defender model was studied in [7], where the attacker aims to maximize the loss while the defender aims to minimize the loss. Further, the trilevel defender-attacker-defender model was investigated in [8]-[10], where the defender develops the optimal defense strategy considering the attacker's optimal strategy in which the power re-dispatch is accounted for. These models consider the interaction between the attacker and defender, but other factors such as the load demand uncertainty and wind generation uncertainty are not considered.

Actually, the load demand and wind generation uncertainties can influence the development of the optimal defense strategy as well. It is important to include them to derive a more comprehensive model for the identification of critical elements, which is the purpose of this paper. The major contribution of this paper is that a defender-attacker-nature-operator model is proposed considering the attacks, the load demand uncertainties, and the wind generation uncertainties. This model can offer more flexibilities for the power system defender to derive a robust defense strategy.

The rest of this paper is organized as follows. Section II presents the mathematical model. The solution method is explained in Section III. Case studies are demonstrated in Section IV to validate the proposed method. Section V summarizes this paper and gives the future work.

II. MATHEMATICAL MODELING

The power system is faced with the risk of deliberate man-made attacks, the generators and transmission lines could all be targeted and attacked. While an obvious way is to protect all the possible targets in case of an attack, it is unrealistic and uneconomical, as there can be a huge number of elements in the power system and all of them could be potential targets. Instead, it is important to identify the most critical components

This project is funded by State Grid Corporation of China (SGCC) under project *Research on Spatial-Temporal Multidimensional Coordination of Energy Resources*.

and minimize the loss in the worst-case scenario. A representative example method is the trilevel model [8]-[10], as shown in Fig. 1.

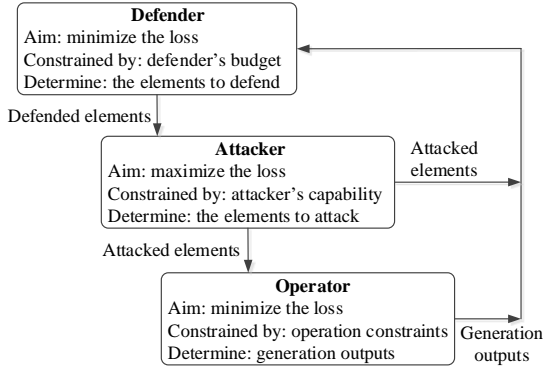

Fig. 1 Conventional trilevel defender-attacker-defender model

There are three agents in this trilevel model. At the top level, the power system defender takes measures to protect some selected components. Specifically, the defender can conduct intensified patrolling, harden the fencing, or deploy alarms to detect and thwart the possible attacks. At the middle level, the attacker chooses some critical components within his/her capacity. At the bottom level, the power system operator re-dispatches the power generation to minimize the loss. These three agents take actions sequentially: the defender takes actions first, followed by the attacker, and the operator acts at last.

The trilevel model assumes the defender can accurately know the loads and renewable generations. Unfortunately, this assumption is challenged by the uncertain nature of the load demands and wind generations, and it is difficult to accurately predict them. If the actual load demands and renewable generations deviate from the predictions, the obtained results based on the fixed predications might not be the optimal. It is necessary and practical to include the load and renewable generation uncertainties in the defender's decision-making.

*A. Load Demand Uncertainty Modeling*

The uncertain load demand $p_i^d$ is expressed as follows.

$$p_i^d = \bar{p}_i^d + \bar{p}_i^{d,+} z_i^{d+} - \bar{p}_i^{d,-} z_i^{d-} \quad \forall i \in I \quad (1)$$

where $\bar{p}_i^d$ is the expected load demand at load point $i$; $\bar{p}_i^{d,+}$ is the maximum upper deviation of load demand $i$; $\bar{p}_i^{d,-}$ is the maximum lower deviation; $z_i^{d+}$ and $z_i^{d-}$ are the load uncertainty factors. The uncertainty factors ($z_i^{d+}$, $z_i^{d-}$) are constrained within a set $\mathbf{Z}^d$, expressed as follows.

$$\mathbf{Z}^d = \begin{cases} 0 \leq z_i^{d+} \leq 1 \ \forall i \in I \\ 0 \leq z_i^{d-} \leq 1 \ \forall i \in I \\ z_i^{d+} + z_i^{d-} \leq 1 \ \forall i \in I \\ \sum_{i \in I}(z_i^{d+} + z_i^{d-}) \leq u^d \end{cases} \quad (2)$$

where $I$ is the set of load demands, and $u^d$ is the budget of load demand uncertainty.

The load demand uncertainty model represented by (1) and (2) suggests that each load can be within a range $[\bar{p}_i^d - \bar{p}_i^{d,-}, \bar{p}_i^d + \bar{p}_i^{d,+}]$, but the realization of the load uncertainty is constrained by $u^d$. A larger value of $u^d$ means more possible realizations. On the contrary, the uncertainty of the load demands is reduced with a smaller value of $u^d$.

*B. Wind Generation Uncertainty Modeling*

Besides the load demands, the wind generation is also stochastic and difficult to be predicted exactly. Assume the set of wind farms is $K$, for each wind farm $k \in K$, the uncertain available wind generation $p_k^{w,avai}$ is modeled as

$$p_k^{w,avai} = \bar{p}_k^w + p_k^{w,+} z_k^{w+} - p_k^{w,-} z_k^{w-} \quad (3)$$

where $\bar{p}_k^w$ is the expected wind generation at wind farm $k$; $p_k^{w,+}$ and $p_k^{w,-}$ are the maximum upper deviation and the maximum lower deviation, respectively; $z_k^{w+}$ and $z_k^{w-}$ are the wind generation uncertainty factors, which are constrained within the set $\mathbf{Z}^w$:

$$\mathbf{Z}^w = \begin{cases} 0 \leq z_k^{w+} \leq 1 \ \forall k \in K \\ 0 \leq z_k^{w-} \leq 1 \ \forall k \in K \\ z_k^{w+} + z_k^{w-} \leq 1 \ \forall k \in K \\ \sum_{k \in K}(z_k^{w+} + z_k^{w-}) \leq u^w \end{cases} \quad (4)$$

The constraints (3)-(4) show that if the uncertainty factor $z_k^{w+}(z_k^{w-})$ is 1, the maximum wind generation upper deviation (lower deviation) is reached. But the realization of all the wind generation uncertainties is limited by the budget of wind generation uncertainty $u^w$.

*C. Proposed Defender-Attacker-Nature-Operator Model*

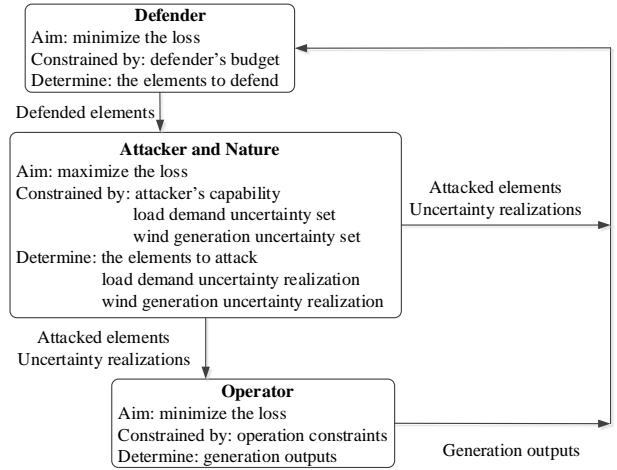

Fig. 2 Proposed defender-attacker-nature-operator model

When the load demand uncertainty and wind generation uncertainty are considered, they can be integrated into the proposed defender-attacker-nature-operator model as shown in Fig. 2. In this model, there are two agents at the middle level, i.e., the attacker and the nature. The nature is an imaginary agent to determine the realization of the load demands and wind generations. These two agents can both influence the power system operation state and the related power re-dispatch performed by the operator. As such, the worst-case scenario is realized when the attacker and the nature collaboratively work to maximize the load loss. Thus, the defender should consider both the intelligent attack strategy and pessimistic realizations of the load demands and wind generations.

The mathematical modeling of this proposed defender-attacker-nature-operator model is described as follows.

$$\min_{w^g,w^f} \max_{\substack{v^g,v^f,z^{w+},\\ z^{w-},z^{d+},z^{d-}}} \min_{\substack{p^g,p^f,\delta,\\ \Delta p^d, p^w}} \sum_{i \in I} \Delta p_i^d \quad (5)$$

$$\sum_{l \in L} c_l^{D,f} w_l^f + \sum_{j \in J} c_j^{D,g} w_j^g \leq r^D \quad (6)$$

$$\sum_{l \in L} c_l^{A,f} (1 - v_l^f) + \sum_{j \in J} c_j^{A,g} (1 - v_j^g) \leq r^A \quad (7)$$

$$p_l^f = \left(w_l^f + v_l^f - w_l^f v_l^f\right) \frac{\sum_{n \in N} A_{nl} \delta_n}{x_l} \quad \forall l \; (\mu_l) \quad (8)$$

$$\sum_{j \in J_n} (w_j^g + v_j^g - w_j^g v_j^g) p_j^g - \sum_{l \in L} A_{nl} p_l^f + \sum_{k \in K_n} p_k^w + \Delta p_{i \in I_n}^d = p_i^d \quad \forall n \; (\lambda_n) \quad (9)$$

$$-\overline{p}_l^f \leq p_l^f \leq \overline{p}_l^f \quad \forall l \; (\underline{\phi}_l, \overline{\phi}_l) \quad (10)$$

$$0 \leq p_j^g \leq \overline{p}_j^g \quad \forall j \; (\overline{\gamma}_j) \quad (11)$$

$$0 \leq p_k^w \leq p_k^{w,avai} \quad \forall k \; (\overline{\beta}_k) \quad (12)$$

$$0 \leq \Delta p_i^d \leq p_i^d \quad \forall i \; (\overline{\alpha}_i) \quad (13)$$

The constants, variable and parameters in (5)-(13) are explained as follows. $w^g$ and $w^f$ are binary variables indicating the defense decisions regarding the generators and branches, respectively. When $w^g$ ($w^f$) is 1, the generator (branch) is protected. $v^g$ and $v^f$ are binary variables indicating the attack decisions regarding the generators and branches, respectively. When $v^g$ is 0, it means the corresponding generator is attacked; otherwise, the corresponding generator is not attacked. And this is the same for $v^f$. $p^g$ means the active power outputs of the generators. $p^f$ refers to the power flow of the branches. $\delta$ indicates the bus voltage angles. $\Delta p^d$ means the load curtailments. $p^w$ refers to the dispatched wind power outputs. $c_l^{D,f}$ and $c_j^{D,g}$ are the costs required to defend branch $l$ and generator $j$, respectively. $c_l^{A,f}$ and $c_j^{A,g}$ are the costs required to attack branch $l$ and generator $j$, respectively. $r^D$ is the defender's capacity, and $r^A$ is the total budget the attacker has. $L$ is the set of branches; $J$ is the set of generators; $N$ is the set of buses. $A_{nl}$ indicates the power flow direction, and it is 1 if the power flow on bus $l$ is defined from bus $n$; it is -1 if the power flow on bus $l$ is to bus $n$; otherwise, it is 0. $x_l$ is the impedance of branch $l$. $J_n$ is the set of generators on bus $n$; $K_n$ is the set of wind generations on bus $n$. $\overline{p}_l^f$ is the maximum power flow capacity of branch $l$; $\overline{p}_j^g$ is the maximum generation capacity of generator $j$.

The objective function is shown in (5) and it is a trilevel optimization problem. Constraint (6) represents the limitation of the defender's budget. Constraint (7) indicate the limitation of the attacker's capacity. Constraints (8)-(13) give the power-dispatch strategy based on the optimal power flow (OPF) analysis. Equation (8) calculates the power flows considering the attack and defense actions. Equation (9) ensures that the incoming and outgoing powers are balanced at each bus. It shows that a branch or a generator will be out-of-service only when it is not protected but attacked. Constraints (10)-(13) indicate the limitations of the power flows, conventional generators, wind power generations, and load curtailments.

### III. SOLUTION METHOD

It is crucial to develop efficient methods to solve the optimization problem represented by (1)-(13). In this paper, the trilevel problem is decomposed into a master problem (MP) and a sub-problem (SP) using the C&CG algorithm [11]. The SP involves the middle level and the bottom level when the decision variables of the top level are given. The MP involves the top level and the bottom level when the decision variables of the middle level are given. The MP and SP are presented in detail as follows.

#### A. Sub-Problem

The SP is for the attacker and nature to make decisions given the defender's defense strategy. Specifically, the attacker determines which elements to attack with the capacity; the nature determines the realization of the load demands and wind generation within the uncertainty sets. The objective function of SP is shown in (14), which is a bilevel problem.

$$\max_{\substack{v^g,v^f,z^{w+},\\ z^{w-},z^{d+},z^{d-}}} \min_{\substack{p^g,p^f,\delta,\\ \Delta p^d, p^w}} \sum_{i \in I} \Delta p_i^d \quad (14)$$

The related constraints are (1)-(4), (7)-(13), in which the values of $\hat{w}^f$ and $\hat{w}^g$ are given. It is noted here that the symbol $\hat{x}$ means the given value of the variable $x$. A widely adopted solution for the bilevel problem is to transform it into an equivalent single-level problem. Also, as shown in [12], the optimal solutions of the uncertain load demands and wind generations will be on the extreme points of their uncertainty sets. Thus, $z_i^{d+}$, $z_i^{d-}$, $z_k^{w+}$ and $z_k^{w-}$ can shrink to binary variables and the issue of multiplication of the continuous variables can be avoided, which can alleviate the computational complexity. The transformed single-level problem is represented in (15)-(30).

$$\eta = \max_{\{v_j^g, v_l^f, z_i^{d+}, z_i^{d-}, z_k^{w+}, z_k^{w-}, \mu_l, \lambda_n, \underline{\phi}_l, \overline{\phi}_l, \overline{\gamma}_j, \overline{\beta}_k, \overline{\alpha}_i\}} \{\sum_{j \in J} \overline{\gamma}_j \; \overline{P}_j^g +$$

$$\sum_{k \in K} \overline{\beta}_k \; (\overline{p}_k^w + p_k^{w,+} z_k^{w+} - p_k^{w,-} z_k^{w-}) + \sum_{l \in L}(\overline{\phi}_l - \underline{\phi}_l) \overline{P}_l^f +$$

$$\sum_{i \in I}(\lambda_{n|i \in I_n} + \overline{\alpha}_i)(\overline{p}_i^d + \overline{p}_i^{d,+} z_i^{d+} - \overline{p}_i^{d,-} z_i^{d-})\} \quad (15)$$

$$\sum_{l \in L} c_l^{A,f}(1 - v_l^f) + \sum_{j \in J} c_j^{A,g}(1 - v_j^g) \leq r^A \quad (16)$$

$$z_i^{d+} + z_i^{d-} \leq 1 \quad (17)$$

$$z_k^{w+} + z_k^{w-} \leq 1 \quad (18)$$

$$\sum_{i \in I}(z_i^{d+} + z_i^{d-}) \leq u^d \quad (19)$$

$$\sum_{k \in K}(z_k^{w+} + z_k^{w-}) \leq u^w \quad (20)$$

$$\sum_{n \in N} A_{nl} \frac{\mu_l(\hat{w}_l^f + v_l^f - \hat{w}_l^f v_l^f)}{x_l} = 0 \quad \forall n \in N \quad (21)$$

$$\lambda_{n|i \in I_n} + \overline{\alpha}_i \leq 1 \quad \forall i \in I \quad (22)$$

$$(\hat{w}_j^g + v_j^g - \hat{w}_j^g v_j^g)\lambda_{n|j \in J_n} + \overline{\gamma}_j \leq 0 \quad \forall j \in J \quad (23)$$

$$\lambda_{n|k \in K_n} + \overline{\beta}_k \leq 0 \quad \forall k \in K \quad (24)$$

$$\mu_l - \sum_{n \in N} A_{nl} \lambda_n + \underline{\phi}_l + \overline{\phi}_l = 0 \quad \forall l \in L \quad (25)$$

$$\overline{\gamma}_j \leq 0 \quad \forall j \in J \quad (26)$$

$$\underline{\phi}_l \geq 0 \quad \forall l \in L \quad (27)$$

$$\overline{\phi}_l \leq 0 \quad \forall l \in L \quad (28)$$

$$\overline{\alpha}_i \leq 0 \quad \forall i \in I \quad (29)$$

$$\overline{\beta}_k \leq 0 \quad \forall k \in K \quad (30)$$

where $\eta$ is the equivalent objective function for (14). In (11)-(30), $\mu_l, \underline{\phi}_l, \overline{\phi}_l, \overline{\gamma}_j, \overline{\beta}_k$ and $\overline{\alpha}_i$ are the duality variables related to constraints (8)-(13). There are nonlinear terms (the product of a binary variable and a continuous variable) in (15), (21) and (23), the big-M method is adopted to linearize them.

### B. Master Problem

Given the middle level decisions $\hat{v}^{f,m}$, $\hat{v}^{g,m}$, $\hat{p}^{d,m}$, $\hat{z}^{d+,m}$, $\hat{z}^{d-,m}$, $\hat{p}^{w,avai,m}$, $\hat{z}^{w+,m}$ and $\hat{z}^{w-,m}$ for each iteration $m = 1, \cdots, m^{iter}$, the MP at iteration $m^{iter}$ is shown by (31)-(41).

$$\min_{w^g, w^f, p^{g,m}, p^{f,m}, \delta^m, \Delta p^{d,m}, p^{w,m}} \xi \quad (31)$$

$$\xi \geq \sum_{i \in I} \Delta p_i^{d,m} \quad \forall m = 1, \cdots, m^{iter} \quad (32)$$

$$\sum_{l \in L} c_l^{D,f} w_l^f + \sum_{l \in L} c_j^{D,g} w_j^g \leq r^D \quad (33)$$

$$\hat{p}_i^{d,m} = \overline{p}_i^d + \overline{p}_i^{d,+} \hat{z}_i^{d+,m} - \overline{p}_i^{d,-} \hat{z}_i^{d-,m} \quad \forall i \in I \quad (34)$$

$$\hat{p}_k^{w,avai,m} = \overline{p}_k^w + p_k^{w,+} \hat{z}_k^{w+,m} - p_k^{w,-} \hat{z}_k^{w-,m} \quad \forall k \in K \quad (35)$$

$$p_l^{f,m} = \left(w_l^f + \hat{v}_l^{f,m} - w_l^f \hat{v}_l^{f,m}\right) \frac{\sum_{n \in N} A_{nl} \delta_n^m}{x_l} \quad \forall l \in L \quad (36)$$

$$\sum_{j \in J_n} \left(w_j^g + \hat{v}_j^{g,m} - w_j^g \hat{v}_j^{g,m}\right) p_j^{g,m} - \sum_{l \in L} A_{nl} p_l^{f,m} + \sum_{k \in K_n} p_k^{w,m} + \Delta p_{i \in I_n}^{d,m} = \hat{p}_i^{d,m} \quad \forall n \in N \quad (37)$$

$$-\overline{p}_l^f \leq p_l^{f,m} \leq \overline{p}_l^f \quad \forall l \in L \quad (38)$$

$$0 \leq p_j^{g,m} \leq \overline{p}_j^g \quad \forall j \in J \quad (39)$$

$$0 \leq p_k^{w,m} \leq \hat{p}_k^{w,avai,m} \quad \forall k \in K \quad (40)$$

$$0 \leq \Delta p_i^{d,m} \leq \hat{p}_i^{d,m} \quad \forall i \in I \quad (41)$$

For conciseness, the iteration $\forall m = 1, \cdots, m^{iter}$ is not shown in (34)-(41). The MP is a mixed integer linear programming (MILP) problem, which can be solved by commercial tools like CPLEX solver.

### C. C&CG Algorithm

Based on the MP and SP, the proposed defender-attacker-nature-operator model is solved using C&CG algorithm [11] illustrated in Fig. 3.

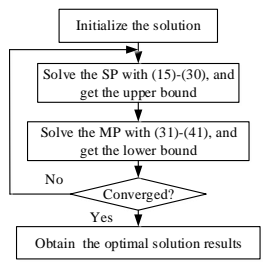

Fig. 3 C&CG algorithm implementation

## IV. CASE STUDIES

The case studies are carried out to validate the proposed model and solution method based on a modified RTS79 system. The original RTS79 system [13] has 24 buses, 2850 MW load and 32 generators. In our case studies, the original generators 3, 14, and 31 are removed, and 3 wind farms are added to buses 1, 13, and 23. The case studies are performed using MATLAB and CPLEX solver.

### A. Base Case

As a case study, assume the expected wind generation at those three wind farms are 160 MW, 150 MW and 120 MW, respectively; the maximum upper deviations and the maximum lower deviations of the wind generations are all 20% of the expected wind generations. Further, the expected values of the load demand predication are as in [12], and the maximum upper deviations and the maximum lower deviations of the load demands are all 30 MW. The budget of load demand uncertainty $u^d$ is 5, and the budget of wind generation uncertainty $u^w$ is 3. $c_l^{D,f}$, $c_j^{D,g}$, $c_l^{A,f}$ and $c_j^{A,g}$ are 1. And the defender's budget is 3 and the attacker's capability is 3.

The simulation is conducted using a laptop with 8 GB RAM, and it takes 419 seconds to finish the simulation. The convergence of the C&CG method is reached after six iterations as shown in Fig. 4.

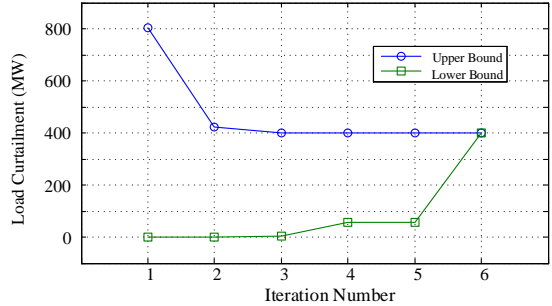

Fig. 4 Convergence of the C&CG method

The simulation results are shown in Table I, also the realizations of the uncertainties are shown in Fig. 5. As a comparison to demonstrate the coordination between the attacker and nature, a case study is conducted when the attack is considered but neglecting the load demand and wind generation uncertainties, the load loss will be reduced to 210 MW and the protected/defended elements can change, as shown in Table I. Also, if the attacks are not considered but the load demand are considered, the load curtailment is 0 MW. These comparisons show that the worst-case load loss can only be realized with the coordination between the attacks and the uncertainties of load demands and wind generations.

Table I. Defended and attacked components

| Variables | Consider uncertainties and attacks | Only consider attacks |
|---|---|---|
| Load loss (MW) | 399 | 210 |
| Defended lines | 25 | 37 |
| Defended generators | 20, 21 | 20, 21 |
| Attacked lines | N/A | 25, 26, 28 |
| Attacked generators | 11, 12, 29 | N/A |

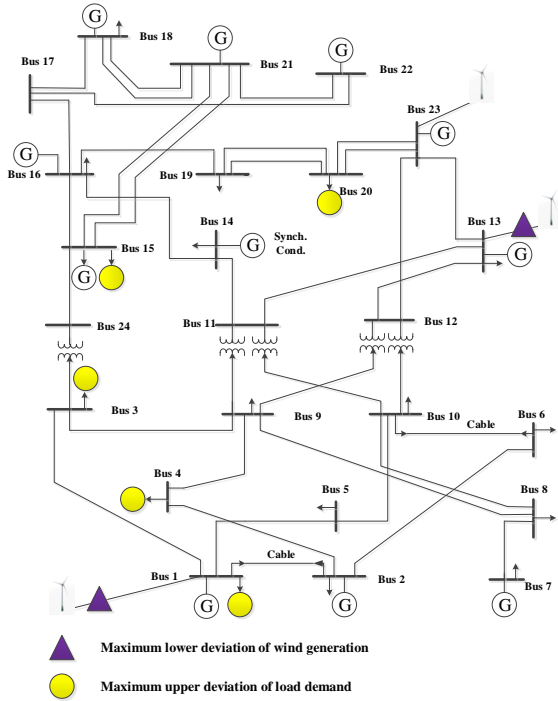

Fig. 5 Realization of the uncertainties in the worst-case scenario

### B. Sensitivity Analysis of Load Demand Uncertainty

The influence of the load demand uncertainty is shown when the maximum upper/lower deviations of the load demands vary while all other parameters are the same as in the base case. The simulation results are shown in Table II. It can be seen that with the increase of deviations, the load loss in the worst-case scenario increases.

Table II. Sensitivity analysis of load demand uncertainty

| Maximum deviations (MW) | Load loss (MW) |
| --- | --- |
| 0 | 272 |
| 10 | 322 |
| 20 | 349 |
| 30 | 399 |
| 40 | 427 |

### C. Sensitivity Analysis of Defense Budget

Table III. Sensitivity analysis of defense budget

| | Defense Budget | | | | | |
| --- | --- | --- | --- | --- | --- | --- |
| | 0 | 1 | 2 | 3 | 4 | 5 |
| Load loss (MW) | 805 | 602 | 422 | 399 | 375 | 291 |
| Defended lines | N/A | N/A | N/A | 25 | 25 | 28, 29 |
| Defended generators | N/A | 21 | 20, 21 | 20, 21 | 20, 21, 29 | 20, 21, 29 |
| Attacked lines | N/A | N/A | N/A | N/A | 23, 27, 29 | 21, 22, 23 |
| Attacked generators | 20, 21, 29 | 12, 20, 29 | 25, 26, 28 | 11, 12, 29 | N/A | N/A |

Table III shows the results with different defense budgets. The load loss, defended lines/generators, and attacked lines/generators for each defense budget are provided. It demonstrates that with the increase of defense budget, the load loss will decrease, which can provide valuable information about the budget needed to maintain the power system risk below a certain level. For example, if the power system defender wants to ensure that the load loss is below 400 MW, the minimum budget required is 3.

## V. CONCLUSIONS AND FUTURE WORK

This paper proposes a defender-attacker-nature-operator model to identify the critical components that should be defended in case of a potential attack, and the load demand uncertainties and renewable generation uncertainties are included in this model. The model is a trilevel optimization problem with the defender at the top level, the attacker and nature at the middle level, and the operator at the bottom level. This problem is solved by the C&CG method, which is a widely used robust optimization solution approach. Case studies and sensitivity analyses are performed based on a modified IEEE RTS79 system, which demonstrates that the proposed model is effective in identifying the key elements for power system defense.

In the future work, more uncertainties related to other components will be included, such as the solar generations, energy storage, electric vehicles, etc. Also, other methods will be explored to further decrease the computational time.